\documentclass[11pt]{amsart}

\usepackage{srcltx}  

\usepackage{geometry}                
\usepackage{graphicx}

\usepackage{amsmath}
\usepackage{amsfonts,amssymb,amsthm,amscd,latexsym,euscript}

\usepackage[all]{xy}
\SelectTips{cm}{10}

\newcommand{\op}{{\ensuremath{\textup{op}}}}

\newcommand{\dgrm}[1]{\ensuremath{\smash{\underset{\widetilde{\hphantom{#1}}}{#1}} \mathstrut}}

\newcommand{\orbit}[1]{\ensuremath{\overset {#1} {\underset \ast \downarrow}}}

\newcommand{\dg}[2]{\ensuremath{\overset{#1} {\underset{#2} \downarrow}}}

\newcommand {\cofib} {\ensuremath{\hookrightarrow}}

\newcommand {\fibr} {\ensuremath{\twoheadrightarrow}}

\newcommand {\trivcofib} {\ensuremath{\tilde\hookrightarrow}}

\newcommand {\we} {\ensuremath{\tilde\rightarrow}}

\DeclareMathOperator{\hocolim}{\textup{hocolim}}
\DeclareMathOperator{\holim}{\textup{holim}}

\newcommand{\cal}[1]{\ensuremath{\mathcal #1}}

\newtheorem {theorem1}{Theorem}[section]
\newtheorem {theorem}[theorem1]{Theorem}

\newtheorem {proposition}[theorem1]{Proposition}
\newtheorem {lemma}[theorem1]{Lemma}
\theoremstyle{definition}
\newtheorem {definition}[theorem1]{Definition}
\newtheorem {example}[theorem1]{Example}
\theoremstyle{remark}
\newtheorem {remark}[theorem1]{Remark}

\newcommand{\cat}[1]{\ensuremath{\EuScript #1}}

\newcommand{\colim}{\ensuremath{\mathop{\textup{colim}}}}

\newcommand{\rarrow}{\rightarrow}
\newcommand{\Id}{\ensuremath{\textup{Id}}}

\newcommand{\mor}[1]{\ensuremath{\textup{mor}(#1)}}

\def\ev{\ensuremath{\textit{ev}}}

\newcommand{\Vopenkas}{Vop\v enka's }

\renewcommand{\hom}{\ensuremath{{\rm hom}}}

\newcommand{\DD}{\cal D}
\newcommand{\sS}{\cal S}

\newtheorem*{VariableNoNum}{{\VariableText}}

\newcommand{\Imap}{\operatorname{Inc}}

\title{Brown representability for space-valued functors}
\author{Boris Chorny}
\address{Department of Mathematics\\
               The University of Haifa at Oranim\\
               Tivon 36006\\
                Israel
}
\email{chorny@math.haifa.ac.il}
\date{\today}                                           

\begin{document}
\begin{abstract}
In this paper we prove two theorems which resemble the classical
cohomological and homological Brown representability theorems. The
main difference is that our results classify contravariant
functors from spaces to spaces up to weak equivalence of functors.

In more detail, we show that every contravariant functor from
spaces to spaces which takes coproducts to products up to homotopy, and
takes homotopy pushouts to homotopy pullbacks is naturally weekly
equivalent to a representable functor.

The second representability theorem states: every contravariant continuous functor from the category of finite simplicial sets to simplicial sets taking homotopy pushouts to homotopy pullbacks is equivalent to the restriction of a representable functor. This theorem may be considered as a contravariant analog of Goodwillie's classification of linear functors \cite{Goo:calc2}.
\end{abstract}

\maketitle

\section{Introduction}
The classical Brown representability theorem \cite{Brown} classifies
contravariant functors from
the homotopy category of pointed connected $CW$-complexes to the category of sets satisfying Milnor's
wedge axiom (W) and Mayer-Vietoris property (MV). 
\begin{itemize}
\item[(W):]  $F(\coprod X_i) = \prod F(X_i)$;
\item[(MV):] $F(D)\to F(B)\times_{F(A)} F(C)$ is surjective for every
  homotopy pushout square
  $\vcenter{\xymatrix@=10pt{A\ar[r]\ar[d]&B\ar[d]\\C\ar[r]&D}}$. 
\end{itemize}

In this paper we address a similar
classification problem, but the functors we classify are the
homotopy functors from spaces to spaces, satisfying (hW) and
(hMV), the higher homotopy versions of (W) and (MV).

\begin{itemize}
\item[(hW):]  $F(\coprod X_i) \simeq \prod F(X_i)$;
\item[(hMV):] $\vcenter{\xymatrix@=10pt{F(D)\ar[r]\ar[d]& F(B)\ar[d]\\
  F(C)\ar[r] & F(A)}} $ is a homotopy pullback for every homotopy
  pushout square
  $\vcenter{\xymatrix@=10pt{A\ar[r]\ar[d]&B\ar[d]\\C\ar[r]&D}}$. 
\end{itemize}

Homotopy functors $F\colon \sS^\op\to\sS$ satisfying (hW) and (hMV) are called
\emph{cohomological} in this paper. Our main result appearing in Theorem~\ref{cohomological-rep} below is that such a functor is naturally weakly equivalent to representable functor.

We should mention right away, that by \emph{spaces} we always mean
simplicial sets in this paper. It is well known that the homotopy category of the unpointed spaces fails to satisfy Brown representability \cite[Prop.~2.1]{Heller}. The enriched framework is more forgiving. Our results are formulated for the unpointed spaces, but they remain valid in the pointed situation too.

Note however, that neither our theorem implies Brown representability,
nor the converse. We assume stronger (higher homotopy) conditions about
the functor, but we also obtain an enriched representability result. 

Nevertheless, our result has a natural predecessor from the Calculus of
homotopy functors. Goodwillie's classification of linear functors
\cite{Goo:calc2} is related to the classical homological Brown representability
in the same way as our representability theorem related to the
cohomological Brown representability. 

The second classification result proved in this paper is ``essentially
equivalent'' to Goodwillie's classification of finitary linear
functors. The difference is that we prove a higher homotopy version
of the homological Brown representability representability in its
contravariant form. Recall \cite{Adams} that every cohomological functor from the category of
compact spectra to abelian groups is a restriction of a representable
functor. We prove a non-stable enriched version of this statement: every
contravariant homotopy functor from finite spaces to spaces satisfying
(hMV) is equivalent to a restriction of a representable functor. Such functors are called \emph{homological}\label{brief-def-homological}.

Although there is no direct implications between our theorem and Goodwillie's classification of linear functors, there is an additional feature that our results share. In both cases every small functor may be approximated by an initial, up to homotopy, representable/linear functor, i.e., both constructions may be viewed as homotopical localizations in some model category of functors. However the collection of all functors from spaces to spaces does not form a locally small category (natural transformations between functors need not form a small set in general). Our remedy to this problem is to consider only \emph{small} functors, i.e., the functors obtained as left Kan extensions of functors defined on a small full subcategory of spaces. 

The method of proof of our results deserves a comment. The Yoneda embedding $Y\colon \cal S \to \cal S^{\op}$ of spaces (=simplicial sets) into the category of small contravariant functors has a left adjoint $Z=\ev_{\ast}$. In this paper we introduce a localization on the category of small contravariant functors such that this pair of adjoint functors becomes a Quillen equivalence, while the local objects are equivalent to the representable functors. In other words, we have a new model for spaces, where every homotopy type is represented by a compact (i.e., finitely presentable) object, the representable functor. Unfortunately our new model of spaces is not class-cofibrantly generated, therefore we can not immediately apply it to the theory of homotopical localizations in spaces. Instead we apply it to the study of representability conditions for small functors. 

We express the property for contravariant functors to satisfy (hW) and (hMV) as a local condition, i.e., such functors become local objects with respect to certain class of maps. We identify this class precisely and argue that the class of object local with respect to those maps is exactly the class of functor equivalent to the representable functors, therefore the localization we constructed is the localization with respect to the class of maps ensuring that the local objects are the cohomological functors.

Therefore, to be equivalent to a representable functor is the same as satisfy the conditions (hW) and (hMV), moreover, every functor has the universal, up to homotopy, approximation by a cohomological functor -- the fibrant replacement in the localized model category, see Remark~\ref{approximation by a cohomology functor} for more details.

We finish our paper with an argument that the new models of spaces, appearing as localizations of class-cofibrantly generated model categories, are not class-cofibrantly generated. This conclusion is quite unexpected, because the localization of a combinatorial model category is always a combinatorial model category (at least under \Vopenkas principle) \cite{CaCho}.
 
\subsection{Acknowledgment} We thank Amnon Neeman for numerous helpful conversation, which led to the results in this paper. We also thank Tom Goodwillie and the anonymous referee for helpful remarks about the early version of this paper.

\section{Model categories of small functors and their localization}\label{model-categories}

The object of study of this paper is homotopy theory of contravariant functors from the category of spaces $\sS$ to~$\sS$. The totality of these functors does not form a category in the usual sense, since the natural transformations between two functors need not form a set in general, but rather a proper class. We choose to treat a sufficiently large subcollection of functors, including all interesting functors and forming a locally small category. The next definition describes elements of a reasonably large subcollection.

\begin{definition}\label{DefineSmall}\label{DefineRepresentable}
  Let $\DD$ be a (not necessarily small) simplicial category. A
  functor $\dgrm X:\DD\to\sS$ is \emph{representable} if there is an object
  $D\in \DD$ such that $\dgrm X$ is naturally equivalent to $R^D$, where
  $R^D(D')=\hom_{\DD}(D, D')$. A functor $\dgrm X\colon \cat C \rarrow
  \cal S$ is called \emph{small} if \dgrm X is a small weighted
  colimit of representables.
\end{definition}

\begin{remark}\label{SmallFunctorCategory}
  G.M.~Kelly \cite{Kelly} calls small functors \emph{accessible} and
  weighted colimits \emph{indexed}. He proves that small functors
  form a simplicial category which is closed under small (weighted)
  colimits \cite[Prop.~5.34]{Kelly}. 
\end{remark}

In order to do homotopy theory we need to work in a category which is
not only cocomplete, but also complete (at least under finite limits).
Fortunately, there is a simple sufficient condition in the situation
of small functors. 

\begin{theorem}
  If $\DD$ is cocomplete, then the category $\sS^{\DD}$ of small
  functors $\DD\to\sS$ is complete.
\end{theorem}

The main technical tool used in the prove of the classification
theorem is the theory of homotopy localizations. More specifically, we
apply certain homotopy localizations in the category of small
contravariant functors  $\cal S^{\cal S^{\op}}$, or in a Quillen
equivalent model category of maps of spaces with the equivariant model
structure \cite{Farjoun, Chorny-Dwyer}. 

Let us briefly recall the definitions and basic properties of the involved model categories. The projective model structure on the small contravariant functors was constructed in \cite{Chorny-Dwyer}. The weak equivalences and fibrations in this model category are objectwise. This model structure is generated by the classes of generating cofibrations and generating trivial cofibrations 
\begin{align*}
I &= \{R_A\otimes \partial \Delta^n \hookrightarrow R_{A}\otimes \Delta^n | A\in \cal S, n\geq 0\}, \\
J &= \{R_{A}\otimes \Lambda^n_k \hookrightarrow R_{A}\otimes \Delta^n | A\in \cal S, n\geq k\geq 0\}.
\end{align*}
The classes $I$ and $J$ satisfy the conditions of the generalized small object argument \cite{pro-spaces}, therefore we refer to this model category as \emph{class-cofibrantly generated}, see \cite[Definition~1.3]{pro-spaces} for the detailed definition and discussion. Note that the representable functors are cofibrant objects and the rest of cofibrant objects are obtained as retracts if $I$-cellular objects. 

Another example of a class-cofibrantly generated model category is given by the equivariant model structure on the maps of spaces $\cal S^{[2]}_\text{eq}$. The central concept of the equivariant homotopy theory is the category of \emph{orbits}. In the category of maps of spaces the subcategory of orbits $\cal O_{[2]}$ is the full subcategory of $\cal S^{[2]}$ consisting of diagrams of the form $\dgrm T = \left(\orbit X\right)$, $X\in\cal S$. Motivation of this terminology and further generalization of the concept of orbit can be found in \cite{DZ}. Equivariant homotopy and homology theories were developed in \cite{DF}. The theory of equivariant homotopical localizations was  introduced in \cite{PhDI}. Weak equivalences and fibrations in the equivariant model category are determined by the following rule: a map $f\colon \dgrm X\to \dgrm Y$ is a weak equivalence or a fibration if for every $\dgrm T\in \cal O_{[2]}$ the induced map of spaces $\hom(\dgrm T, f)\colon \hom(\dgrm T, \dgrm X)\to \hom(\dgrm T, \dgrm Y)$. 

The categories of maps of spaces and small contravariant functors are related by the functor $\cal O\colon \cal S^{[2]}\to \cal S^{\cal S^\op}$, called the \emph{orbit-point functor} (generalizing the fixed-point functor from the equivariant homotopy theory with respect to a group action), which is defined by the formula $(\dgrm X)^{\cal O}(Y) = \hom\left( \orbit Y, \dgrm X \right)$, for all $Y\in \cal S$. Orbit-point functor has a left adjoint called the \emph{realization functor} $|-|_{[2]}\colon \cal S^{\cal S^\op} \to  \cal S^{[2]}$. The main result of \cite{Chorny-Dwyer} is that this pair of functors is a Quillen equivalence.

Before proving the main classification result, we suggest the
following alternative characterization of functors satisfying (hW) and
(hMV) as local objects with respect to some class of maps. 

\subsection*{Homotopy functors as local objects} By definition every
cohomological functor $F$ is a homotopy functor, i.e.,
$F(f)\colon F(B)\to F(A)$ is a weak equivalence for every weak
equivalence $f\colon A\to B$. Denote by $\cal F_{1}$ the class of maps
between representable functors induced by weak equivalences: 
\[
\cal F_{1} = \{ f^{\ast}\colon R_A \to R_{B} | f\colon A\to B \,
\text{is a w.e.}\}, 
\]
where $R_{A}$ denotes the representable functor $R_{A} = \cal S(-,A)$.

Yoneda's lemma implies that $\cal F_{1}$-local functors are precisely
the fibrant homotopy functors. 

\subsection*{Cohomology functors as local objects} Given a homotopy
functor $F$, it suffices  to demand two additional properties for the
functor $F$ to be cohomological: $F$ must convert coproducts to
products up to homotopy and it also must convert homotopy pushouts to
homotopy pullbacks. Yoneda's lemma and the standard commutation rules
of various (ho)(co)limits with $\hom(-,-)$ implies that both
properties are local with respect to the following classes of maps: 
\[
\cal F_{2} = \left\{\left. \coprod R_{X_{i}}\to R_{\coprod X_{i}}
\right | \forall \{X_{i}\}_{i\in I}\in \cal S^{I} \right\} 
\]
and
\[
\cal F_{3} = \left\{\left. \hocolim\left( 
                   \vcenter{\xymatrix@=10pt{R_{A}\ar[r]\ar[d] & R_{C}\\
                                                R_{B}}}\right)
\longrightarrow R_D            
                               %
                                                \right |
		\, \vcenter{\xymatrix@=10pt{A\ar[r]\ar[d]& C\ar[d]\\
						 							 B\ar[r] & D}} \text{-- homotopy pushout in } \cal S \right\}.
\]

Objects which are local with respect to $\cal F = \cal F_1 \cup \cal
F_2 \cup \cal F_3$ are precisely the fibrant homotopy functors.

\begin{lemma}\label{homotopy_pullback}
Any functor $F\colon \cal S^{\op}\to \cal S$ satisfying (hMV) is a homotopy functor, i.e., for any weak equivalence $f\colon A\to B$, the map $F(f)\colon F(B)\to F(A)$ is a weak equivalence.
\end{lemma}
\begin{proof}
Given a weak equivalence $f\colon A\to B$ the following commutative square is a
homotopy pushout:
\[
\xymatrix{A\ar@{=}[r]\ar@{=}[d]& A\ar[d]^f\\
 A\ar[r]_f & B.}
\]
Applying $F$ we obtain:
\[
\xymatrix{
F(B)\ar[r]^{F(f)}\ar[d]_{F(f)}& F(A)\ar@{=}[d]\\
F(A)\ar@{=}[r] & F(A).}
\]
The later square is a homotopy pullback iff $F(f)$ is a weak equivalence. Therefore, any functor satisfying (hMV) is automatically a homotopy
functor.
\end{proof}

We conclude that it suffices to invert $\cal F = \cal F_2 \cup \cal F_3$ .

\begin{remark}\label{reduced}
The indexing category $I$ used to describe $\cal F_{2}$ is a completely arbitrary small discrete category. In particular $I$ can be empty. This implies that the map $\emptyset\to R_{\emptyset}$ is in $\cal F_{2}$. In other words, if $F$ is a cohomological functor, then $F(\emptyset) = \ast$. This property is analogous to the requirement that every linear functor is reduced in homotopy calculus.
\end{remark}

\begin{remark}\label{reduced-homological}
Since homological functors (see a brief explanation on p.~\pageref{brief-def-homological} or an official Definition~\ref{def-homological}) are defined on the category of finite simplicial sets, we need to adjust the definition of $\cal F_{3}$. 
\[
\cal F'_{3} = \left\{\left. \hocolim\left( 
                   \vcenter{\xymatrix@=10pt{R_{A}\ar[r]\ar[d] & R_{C}\\
                                                R_{B}}}\right)
\longrightarrow R_D            
                               %
                                                \right |
		\, \vcenter{\xymatrix@=10pt{A\ar[r]\ar[d]& C\ar[d]\\
						 							 B\ar[r] & D}} \text{-- homotopy pushout,}\, A,B,C,D\in \cal S_\text{fin} \right\}.
\]
Then the reduced homological functors in $\cal S^{\cal S^\op_\text{fin}}$ (with the projective model structure) are precisely the functors which are local with respect to $\cal F'=\cal F'_{3}\cup \{\emptyset\to R_{\emptyset}\}$
\end{remark}

\subsection{Localization}\label{Q-localization} Representing the class of cohomology  functors as local objects does not contribute much to their understanding. Our next goal is to make sure that there exists a localization of the model structure with respect to $\cal F$ and the class of objects we a willing to classify will be represented, up to homotopy, by the elements of the homotopy category of the localized model category. After we achieve this, we have a chance to find a simpler model category, Quillen equivalent to the localized model category, hence classifying the objects of the homotopy category. In addition the localization approach to the classification problem supplies us with an approximation tool, namely the fibrant replacement in the localized category, so that every functor may be turned into a cohomological functor in a functorial way and such approximation is initial with respect to maps into other cohomological functors.

Localization procedure is not always a routine. For example, the existence of localization of spaces with respect to the class of cohomological equivalences is still an open problem (assuming \Vopenkas principle in addition to the standard axioms this question was positively settled \cite{CSS}). In our situation no currently existing general localization machine may be immediately applied, since $\cal F$ is a proper class of maps and the category of small functors is not cofibrantly generated. We will implement an ad hoc approach to this localization problem. Namely, relying on the intuition stemming out of the classical Brown representability we assume that the localized model category will be equivalent to the category of spaces, construct such localization disregarding $\cal F$, and afterwards prove that this localization is precisely the localization with respect to $\cal F$.  

The basic idea behind this localization is to turn the adjunction $\ev_{*}\colon \cal S^{\cal S^{\op}}\rightleftarrows \cal S:\! Y$ into a Quillen equivalence (to see that this is indeed an adjunction note that $\ev_{*}(F)=F\star \Id_{\cal S}$). For this purpose we will use the derived version of the unit of this adjunction: $F\to Y\ev_{*}F$. We need to turn $q=Y\ev_{*}\colon \cal S^{\cal S^{\op}}\to\cal S^{\cal S^{\op}}$ into a homotopy functor.

Since $\ev_{*}$ is a homotopy functor in the projective model structure and $Y$ preserves weak equivalences of fibrant simplicial sets, the derived version of $q$ may be chosen to be the composition $Q=Y\widehat{\ev_*}$, where $\widehat{(-)}$ is a functorial fibrant replacement in simplicial sets. $Q$ is equipped with a coaugmentation $\eta\colon\Id\to Q$, defined as a composition of the unit of adjunction with the application of $Y$ on the natural map of simplicial sets $\ev_{*}(F)\to\widehat{\ev_*(F)}$.

The category $\cal S^{[2]}_\text{eq}$ is related to the category of contravariant functors by the Quillen equivalence \cite{Chorny-Dwyer}:
\begin{equation}
|-|_{[2]}\colon \cal S^{\cal S^\op}\rightleftarrows \cal S^{[2]}_\text{eq}:\!(-)^{\cal O}. \label{Quillen-equiv}
\end{equation}
We would like to localise simultaneously the model category $S^{[2]}_\text{eq}$, so that the adjunction (\ref{Quillen-equiv}) would remain Quillen equivalence.

In order to construct the required localization of $\cal S^{[2]}_\text{eq}$ we will take the derived version of the unit of the adjunction
\begin{equation}
L\colon \cal S^2_\text{eq} \leftrightarrows \cal S :\!R,
	\label{adjunction}
\end{equation}
where $L\left(\overset A {\underset B \downarrow}\right) = A$ and $R(A) = \orbit A$. We define $Q'\left(\overset A {\underset B \downarrow}\right) = \orbit {\hat A}$, and notice that the unit of the adjunction (\ref{adjunction}), composed with the application of $R$ on the natural map $L\left(\overset A {\underset B \downarrow}\right) \to \widehat{L\left(\overset A {\underset B \downarrow}\right)}$, provides $Q'$ with a coaugmentation $\eta'\colon \Id\to Q'$.

It turns out that the localization of the model category $\cal S^{[2]}_\text{eq}$ with respect to $Q'$  is precisely the localization of  $\cal S^{[2]}_\text{eq}$ with respect to the class of maps $|\cal F|_{[2]} = |\cal F_1|_{[2]}\cup|\cal F_2|_{[2]}\cup|\cal F_3|_{[2]}$, where 
\begin{eqnarray*}
|\cal F_1|_{[2]} = \left\{\left . \orbit A  \longrightarrow \orbit B \right | A\to B \text{ is a w.e. in } \cal S \right\},\\
|\cal F_2|_{[2]} = \left\{\left . \coprod \orbit{X_i} \longrightarrow \orbit {\coprod X_i} \right | \forall \{X_{i}\}_{i\in I}\in \cal S^{I} \right\},
\end{eqnarray*}
and 
\[
|\cal F_3|_{[2]} = \left\{
\left. 
\hocolim\left( 
                   \vcenter{\xymatrix@=10pt{{\orbit A}\ar[r]\ar[d] & {\orbit C}\\
                                                {\orbit B}
                                               }                                           
                   }
        \right)\longrightarrow 
        \orbit{D
                 }
\right| 
\vcenter{
\xymatrix@=10pt{A\ar[r]\ar[d]& C\ar[d]\\
			    B\ar[r]        & D
	     }
} \text{ is a homotopy pushout in } \cal S
\right\}.
\]

\begin{remark}
The realization functor $|-|_{[2]}$ may be viewed as a coend $\Imap\otimes_{\cal S} - $, where $\Imap\colon \cal S = \cal O_{[2]} \hookrightarrow \cal S^{[2]}$ is the fully-faithful embedding of the subcategory of orbits \cite{Chorny-Dwyer}. Therefore, computing the realization of the representable functors is just the evaluation of $\Imap$ at the representing object, since the dual of the Yoneda lemma applies.
\end{remark}

The main technical achievement of this paper, which is behind the proof of the representability theorem is the following.
\begin{theorem}\label{main-theorem}
There exist localizations of the projective model structure on $\cal S^{\cal S^{\op}}$ with respect to $Q$ and of the equivariant model structure on $\cal S^{[2]}$ with respect to $Q'$, so that all adjunctions in the following triangle become Quillen equivalences.
\[
\xymatrix{
&\cal S\ar@/^/[ddr]^{R} \ar@/^/[ddl]^{Y} \\
\\
\cal S^{\cal S^\op} \ar@/^/[rr]^{|-|_{[2]}} \ar@/^/[uur]^{\ev_{*}} & &\cal S^{[2]} \ar@/^/[ll]^{(-)^{\cal O}}\ar@/^/[uul]^{L}
}
\]
\end{theorem}

\begin{proof}
The existence of localization follows from the Bousfield-Friedlander theorem \cite[A.7]{BF:gamma}. We have to verify that 
\begin{enumerate}
\item $Q$ and $Q'$ preserve weak equivalences;
\item $Q$ and $Q'$ are coaugmented, homotopy idempotent functors;
\item Pull back of a $Q$($Q'$)-equivalence along a $Q$($Q'$)-fibration is a $Q$($Q'$)-equivalence again (the resulting localized category becomes right proper).
\end{enumerate}
$Q$ and $Q'$ are constructed in such a way that conditions (1) and (2) are satisfied. The verification is a routine.

In order to verify (3) notice that a map in $\cal S^{\cal S^{\op}}$ ($\cal S^{[2]}$) is a $Q$($Q'$)-equivalence iff the map induced between the values of the functors in $*\in\cal S$ ($0\in [2]$) is a weak equivalence. Since $\cal S$ is right proper, any pull back of such map along a levelwise fibration will have the same property. Certainly any $Q$($Q'$)-fibration is a levelwise fibration, hence the conditions of Bousfield-Friedlander theorem are satisfied.

It remains to show that the adjunctions in the triangle above became Quillen equivalences. Note that the composition of the right adjoints of the right edge and of the base of the triangle equals to the right adjoint of the left edge $(R(-))^{\cal O} = Y(-)$, so it suffices to verify only that the right edge and the base of the triangle are Quillen equivalences.

The adjunction of the right edge is a Quillen pair, since the left adjoint $L$ preserves cofibrations and trivial cofibrations. It remains to show that $A=L\left(\dg{A}{B}\right)\to X$ is a weak equivalence iff $\left(\dg A B\right) \to R(X)=\left(\orbit X\right)$ is a $Q'$-equivalence, which is clear.

The adjunction in the base of the triangle is a Quillen pair by Dugger's lemma \cite[8.5.4]{Hirschhorn}, since the right adjoint preserves fibrations of fibrant objects (in the category of maps $Q'$-fibrant object are weakly equivalent to orbits,  hence their orbit points are weakly equivalent to representable functors, i.e., $Q$-fibrant in the category of contravariant functors, but fibrations of $Q$-local objects are $Q$-fibrations), and also trivial fibrations (since those do not change under localization). It remains to show that for  all cofibrant $F\in \cal S^{\cal S^{\op}}$ and for all $Q'$-fibrant $\left(\dg A B\right) \,\simeq\, \left(\orbit A\right)$ a map $f\colon F\to\left(\orbit A\right)^{\cal O}=R_{A}$ is a $Q$-equivalence iff the adjoint map $f^{\sharp}\colon |F|_{[2]}\to \left(\orbit A\right)$ is a $Q'$-equivalence. The `only if' direction follows by applying the realization functor on $f$, since $|R_{A}|_{[2]}=\,\left(\orbit A\right)$ and realization preserves weak equivalence of cofibrant objects. The `if' direction follows from computation of the value $F(*)$: using the composition of two left adjoints $L(|F|_{[2]})=\ev_{*}(F)$, we find out that $F(*)$ is equivalent to the domain of $|F|_{[2]}$).
\end{proof}

\begin{remark}
Theorem~\ref{main-theorem} provides us with two model of spaces with the following property: every object is weakly equivalent to an $\aleph_{0}$-small object. This conclusion seems contra-intuitive in view of Hovey's proof that every cofibrant and $\aleph_0$-small, relative to cofibrations, object in a pointed finitely generated model category \cat C is $\aleph_0$-small in $\mathrm{Ho}(\cat C)$ \cite[7.4.3]{Hovey}. However, there is no contradiction with our result, since the localized model categories $\cal S^{[2]}_\text{eq}$ or $\sS^{\sS^{\op}}$ are very far from being finitely generated. 

It is tempting to try to apply these models to the problem of localization of spaces with respect to some proper class of maps, which we could not do before due to set theoretical difficulties (the cardinality of domains and codomains of these maps would not be bounded by any fixed cardinal).  However, there is still an obstacle preventing an immediate application of these models to localization problems in $\cal S$. The  Bousfield-Friedlander localization machinery used to prove Theorem~\ref{main-theorem} does not provide the localized model categories with a class of generating trivial cofibrations that is necessary for construction of new localizations. In fact the new model  categories fail to be (class-)cofibrantly generated, as we will show in Section~\ref{non-cofib}.
\end{remark}

Our next goal is to show that $Q$-localization is precisely the localization with respect to $\cal F$ and $Q'$-localization is precisely the localization with respect to $|\cal F|_{[2]}$.

\section{Technical preliminaries}
Recall that we are going to prove two more theorems in this paper. Theorem~\ref{cohomological-rep} classifies cohomological functors and Theorem~\ref{homological-rep} classifies homological functors. However the technicalities behind the proofs are very similar. Therefore, while we are heading towards the proof of Theorem~\ref{cohomological-rep} first, we indicate little adjustments required to adapt the argument for the proof of Theorem~\ref{homological-rep}.

The $Q$-local objects are precisely the functors (levelwise) weakly equivalent to the representable functors  $R_A$ with $A$ fibrant. We need to show that every object in $\cal S^{\cal S^{\op}}$ is $\cal F$-local equivalent to a representable functor.

Every small contravariant functor may be approximated by an $I$-cellular diagram, up to a (levelwise) weak equivalence \cite{Chorny-Dwyer}, where 
\[
I = \left\{\left.
{\overset {\partial \Delta^n} {\underset {\Delta^n} \downarrow}} \otimes R_{A} \,
\right| \,  A\in \cal S \right\}.
\]
Therefore, it suffices to show that every $I$-cellular diagram is $\cal F$-equivalent to a representable functor. We are going to prove it by cellular induction, but we precede the proof with the following lemma, which says that the basic building blocks of cellular complexes are $\cal F$-equivalent to representable functors.

\begin{lemma}\label{sphere}
For every $A\in \cal S$, $n\ge 0$, there exists $A'\in \cal S$ such that $\partial\Delta^{n}\otimes R_{A}\overset {\cal F} \simeq R_{A'}$.
\end{lemma}
\begin{proof}
We will prove the statement with $A'\simeq \partial\Delta^n\otimes A$. The proof is by induction on $n$. For $n = 0$ we have $\partial\Delta^{0}\otimes R_{A} = \emptyset \otimes R_{A} = \emptyset \overset {\cal F} \simeq R_{\emptyset} = R_{\partial\Delta^0\otimes A}$, since the map $\emptyset \to R_{\emptyset}$ is in $\cal F$ by Remark~\ref{reduced}. Alternatively, if one is willing to exclude $\cal F_2$ from $\cal F$, then for the base of induction it suffices to assume that the cohomology functor $F$ is reduced, i.e., $F(\emptyset) = \ast$; cf. Remark~\ref{reduced-homological}. In other words the basis for induction holds for $\cal F'$ equivalences as well.

Suppose the statement is true for $n$, i.e., $\partial\Delta^{n}\otimes R_{A} \overset {\cal F}\simeq R_{\partial\Delta^{n}\otimes A}$; we need to show it for $n+1$.
\begin{multline*}
\partial\Delta^{n+1}\otimes R_{A} \simeq 
\colim 
\left(
 \vcenter{
  \xymatrix{\partial \Delta^n\otimes R_{A} \ar@{^(->}[r]\ar@{^(->}[d]& \Delta^n\otimes R_{A}\\
  					\Delta^n\otimes R_{A}
  }
 } 
\right)\overset {\cal F} \simeq
\hocolim 
\left(
 \vcenter{
  \xymatrix{R_{\partial \Delta^n\otimes A} \ar[r]\ar[d]& R_{A}\\
  					 R_{A}
  }
 } 
\right)\simeq\\
\hocolim 
\left(
 \vcenter{
  \xymatrix{ R_{\partial \Delta^n\otimes A} \ar[r]\ar[d]& R_{\Delta^n\otimes A}\\
  					 R_{\Delta^n\otimes A}
  }
 } 
\right)\overset {\cal F} \simeq 
  R_{
    \colim 
      \left(
        \vcenter{
            \def\objectstyle{\scriptstyle}
            \xymatrix@=12pt{ {\partial \Delta^n\otimes A} \ar@{^(->}[r]\ar@{^(->}[d]& \Delta^n\otimes A\\
  	  			     	  \Delta^n\otimes A
            }
       }
      \right)
   }
\simeq
R_{(\Delta^n\coprod_{\partial\Delta^n}\Delta^n)\otimes A} \simeq R_{\partial\Delta^{n+1}\otimes A},
\end{multline*}
where  the first $\cal F$-equivalence is induced by the $\cal F$-equivalence in the upper left vertex of the diagram (by induction hypothesis) and in the other two vertices we have levelwise weak equivalences.
(If we will map both homotopy pushouts into an arbitrary $\cal F$-local object \dgrm W, we will obtain a levelwise weak equivalence of homotopy pullback squares of spaces). The second $\cal F$-equivalence is induced by the map from $\cal F_3\subset \cal F$ corresponding to the homotopy pushout square:
\[          
            \xymatrix{ {\partial \Delta^n\otimes A} \ar@{^(->}[r]\ar@{^(->}[d]& \Delta^n\otimes A\ar[d]\\
  	  			     	  \Delta^n\otimes A\ar[r] & (\Delta^n\coprod_{\partial\Delta^n}\Delta^n)\otimes A.
            }
\]

The above argument applies for all finite $A$ if we consider $\cal F'$ instead of $\cal F$, as we did not use any equivalences induced by an element of $\cal F_{2}$.
\end{proof}

We will need to use the following standard result 
\begin{lemma}\label{pushout-trick}
The following commutative square is a pushout square
\[
\xymatrix{
A
\ar[r]^{f}
\ar[d]_{g}
		 & B
		   \ar[d]^{g'}\\
C 
\ar[r]_{f'}
		& D
}
\]
if and only if the square
\[
\xymatrix{
A\coprod A
\ar[r]^{\nabla}
\ar[d]_{f\coprod g}
		 & A
		   \ar[d]^{g'f}\\
B\coprod C 
\ar[r]
		& D
}
\]
is a pushout square.
\end{lemma}
\begin{proof}
Represent the two pushout diagrams as the coequalizers:
\[
\xymatrix{
A\coprod A 
\ar[r]<1ex>^<<<<<{\nabla}
\ar[r]<-1ex>_<<<<{f\coprod g} & A\coprod B\coprod C \ar[r] & D
}
\]
and
\[
\xymatrix{
A\coprod A \coprod A \coprod A 
\ar[r]<1ex>^<<<<<{\nabla^2}
\ar[r]<-1ex>_<<<<{(f\coprod g)\coprod(f\coprod g)} & A\coprod B\coprod C \ar[r] & D
}
\]
There exist natural maps in both directions between the coequalizer diagrams, showing that their colimits coincide.
\end{proof}

\begin{lemma}\label{obvious-lemma}
Let $\cat M$ be a class-cofibrantly generated model category, such that the class of generating cofibrations $I$ has $\aleph_0$-small domains with respect to the cofibratons. Then every $I$-cellular complex $X\in \cat M$ may be decomposed into an $\omega$-indexed colimit $X=\colim_{n}X_{n}$ such that for every $n\in \mathbb N$ there is a pushout square
\begin{equation}\label{few-cell-attach}
\xymatrix{
A
\ar[r]
\ar@{^{(}->}[d]_{f} 
  &  X_n
       \ar@{^{(}->}[d]\\
B
\ar[r]
  &  X_{n+1},
}
\end{equation}
where the map $A\cofib B$ is a coproduct of a set of maps from $I$.
\end{lemma}
\begin{proof}
Every $I$-cellular complex $X$ has a decomposition into a colimit indexed by a cardinal $\lambda$:
\[
X = \colim_{a<\lambda}(X_{0,0}\to\cdots\to X_{a,0}\to X_{a+1,0}\to \cdots),
\]
where $X_{0,0}=\emptyset$, $X_{a,0}$ is obtained from $X_{a-1,0}$ by attaching a cell $g\in I$:
\begin{equation}\label{cell-attach}
\xymatrix{
C \ar[r]\ar@{^(->}[d]_{g} & X_{a-1,0}\ar@{^(->}[d]\\
D \ar[r] &  X_{a,0},}
\end{equation}
and  $X_{a,0}=\colim_{b<a}X_{b,0}$ for all limit ordinals $a$.

If $\lambda=\omega$, then we are done, otherwise assume for transfinite induction that for all $\omega \leq b < a < \lambda$
\[
X_{b,0} = \colim_{n\in \mathbb N} X_{n,b}, 
\]
so that for all $n\in \mathbb N$, $X_{n+1,b}$ is obtained from $X_{n,b}$ by attaching a coproduct of a set of maps from $I$.

If $a$ is a successor cardinal, then $X_{a,0}$ is obtained from $X_{a-1,0}$ as a pushout (\ref{cell-attach}). On the other hand, by the inductive assumption, $X_{a-1,0}=\colim_{n<\omega} X_{n,a-1}$. Since the domain $C$ of the cell $g\in I$ is $\aleph_0$-small with respect to cofibrations, the attaching map $C\to X_{a-1,0}$ factors through $X_{k,a-1}$ for some finite $k$.

Let us define $X_{n,a} = X_{n,a-1}$, if $n\leq k$. If  $n=k+1$, then we define first $X'_{k,a}=D\coprod_{C} X_{k+1,a-1}$, and now $X_{k+1,a}=X'_{k,a}\coprod_{X_{k,a-1}}X_{k+1,a-1}$. And if $n> k+1$, then put $X_{n,a}=X_{n-1,a}\coprod_{X_{n-1,a-1}}X_{n, a-1}$. We have $X_{a,0}=\colim_{n<\omega} X_{n,a}$, since in the commutative diagram
\[
\xymatrix{
C
\ar[r]
\ar@{^(->}[d]_{g}
		& X_{k,a-1}
		    \ar[d]\ar[r]
				& X_{k+1,a-1}
				   \ar[d]\ar[r]
						& \ldots 
						   \ar[r]
								& X_{a-1,0}=\colim_{n<\omega}X_{n,a-1}
								   \ar@<-1,5cm>[d]\\
D
\ar[r]
		& X'_{k,a}
		   \ar[r]
				& X_{k+1,a}
				   \ar[r]
						& \ldots
						   \ar[r]
								& X_{a,0}=\colim_{n<\omega}X_{n,a}
}
\]
all squares composing the ladder are pushouts by definition, so is the outer square.

It remains to show that $X_{k+1,a}$ is obtained from $X_{k,a}$ as in a pushout of the form  (\ref{few-cell-attach}). For other values of the first index this is clear. It suffices to show that the square
\[
\xymatrix{
A\coprod C
\ar@{^(->}[d]
\ar[r] & 
		X_{k,a} 
		\ar@{=}[r] 
		\ar@{^(->}[d]&
				X_{k,a-1}\\
B\coprod D 
\ar[r]&
		 X_{k+1,a}
}
\]
is a pushout.
First let us split it into two squares
\begin{equation}\label{two-squares}
\xymatrix{
A\coprod C
\ar@{^(->}[d]
\ar[r] &
	A\coprod X_{k,a-1}
	\ar@{^(->}[d]
	\ar[r] &
		X_{k,a} 
		\ar@{^(->}[d]\\
B\coprod D 
\ar[r]&
	B\coprod X'_{k,a}
	\ar[r] &
		 X_{k+1,a}
}
\end{equation}
and then show that these two squares are pushouts.

The left square is a pushout as a coproduct of two pushout squares. It remains to show that the right square of (\ref{two-squares}) is a pushout. 

Let us start with the following pushout square:
\[
\xymatrix{
A
\ar[r]
\ar@{^(->}[d] & 
	X'_{k,a}
	\ar@{^(->}[d]\\
B
\ar[r] & 
	X_{k+1,a}.
}
\]
Lemma~\ref{pushout-trick} implies that the square
\[
\xymatrix{
A\coprod A 
\ar[d]
\ar[r]&
		B\coprod X'_{k,a}
		\ar[d]\\
A 
\ar[r]& 
		X_{k+1,a}
}
\]
is also a pushout. Now we can split it into two squares again
\[
\xymatrix{
A\coprod A 
\ar[r]
\ar@{^(->}[d] & 
		A\coprod X_{k,a-1} 
		\ar@{^(->}[d] 
		\ar[r]& 
				B\coprod X'_{k,a}
				\ar[d]\\
A 
\ar[r]& 
		X_{k,a}
		\ar[r]& 
				X_{k+1,a},
}
\]
where the right square is exactly the right square of (\ref{two-squares}) and the left square is a pushout by Lemma~\ref{pushout-trick}  since the square
\[
\xymatrix{
A 
\ar@{=}[d]
\ar[r]&
		 X_{k,a-1}
		 \ar@{=}[d]\\
A 
\ar[r]& 
		X_{k,a}
}
\]
is a pushout. Therefore the right square is also a pushout, which is what we needed to show. 
\end{proof}

\begin{proposition}\label{main_prop}
Every $I$-cellular complex $\dgrm X\in \cal S^{\cal S^{\op}}$ is $\cal F$-equivalent to a representable functor $R_A$ for some $A$.
\end{proposition}
\begin{proof}
By definition, every $I$-cellular complex $X$ has a decomposition into a colimit indexed by a cardinal $\lambda$ starting from the initial object and on each stage one element of $I$ is attached. By Lemma~\ref{obvious-lemma}, there is an alternative decomposition of $X$:
\[
X = \colim_{a<\omega}(X_{0}\to\cdots\to X_{a}\to X_{a+1}\to \cdots),
\]
where $X_0=\emptyset$ and $X_{a+1}$ is obtained from $X_{a}$ by attaching a small collection cells:
\[
\xymatrix{
\coprod_{A}(\partial\Delta^{n}\otimes R_A) \ar[r]\ar@{^(->}[d] & X_{a}\ar@{^(->}[d]\\
\coprod_{A}(\Delta^{n}\otimes R_A) \ar[r] &  X_{a+1}.}
\]

Note for the basis of induction, that $X_{1}$ is $\cal F$-equivalent to $R_{\coprod_{A}A}$, since $\Delta^{n}\otimes R_A \simeq R_A$. Assuming, by induction, that $X_{a}$ is $\cal F$-equivalent to a representable functor $R_{C_a}$, we notice, by Lemma~\ref{sphere}, that $\coprod_{A}(\partial\Delta^{n}\otimes R_A)=\partial\Delta^{n}\times \coprod_{A}R_A$ is $\cal F$-equivalent to $R_{\partial\Delta^{n}\otimes \coprod_{A} A}=R_{\coprod_{A} (\partial\Delta^{n}\otimes A)}$, and $\coprod_{A}(\Delta^{n}\otimes R_A) \simeq R_{\coprod_{A}A}$, so all the vertices of the homotopy pushout above are $\cal F$-equivalent to representable functors $R_{A'}$ for some $A'$. We conclude that $X_{a+1}$ is $\cal F$-equivalent to a representable functor $R_{C_{a+1}}$, where $C_{a+1}$ is the homotopy pushout $(\coprod_{A}A\leftarrow \coprod_{A}(\partial \Delta^{n} \otimes A) \to C_a)$, similarly to the argument of Lemma~\ref{sphere}.

We obtain the following countable commutative ladder:
\[
\xymatrix{
X_{0}\ar@{^(->}[r] 
       \ar[d]^{\dir{~}}_{\cal F}&
                                          \cdots \ar@{^(->}[r]   & X_{a}\ar@{^(->}[r] 
                                                                               \ar[d]^{\dir{~}}_{\cal F}& 
                                                                                                                          X_{a+1}\ar@{^(->}[r] 
                                                                                                                          \ar[d]^{\dir{~}}_{\cal F}& \cdots\\
R_{C_0} \ar[r]             & \cdots \ar[r]               & R_{C_{a}}\ar[r]             & R_{C_{a+1}}\ar[r]             & \cdots
}.
\]
Taking homotopy colimit of the upper and the lower rows we find that $X \overset{\cal F}\simeq \hocolim_{a<\omega} R_{C_{a}}$, since if we will map both homotopy colimits into an arbitrary $\cal F$-local functor $W$, we will obtain a weak equivalence between the homotopy inverse limits.

Finally, $\hocolim_{a<\omega} R_{C_{a}} = \hocolim_{a<\omega} \left(R_{C_0}\stackrel{f_0} {\longrightarrow} \cdots\longrightarrow R_{C_a}\stackrel{f_a}{\longrightarrow} R_{C_{a+1}}\stackrel{f_{a+1}}{\longrightarrow} \cdots\right)$ may be represented as a homotopy pushout as follows:
\[
\hocolim_{a<\omega} R_{C_{a}} \simeq
\hocolim\left(
\vcenter{
\xymatrix{
\left(\coprod R_{C_a}\right)\coprod
\left(\coprod R_{C_a}\right)  \ar[r]^(.65){1\coprod f}\ar[d]^\nabla & \coprod R_{C_a}\\
\coprod R_{C_a}
}
}
\right),
\]
where $f = \coprod_{a<\omega}{f_a}$ is the shift map and $\nabla$ is the codiagonal. Observe that the homotopy pushout above is weakly equivalent to the infinite telescope construction.

All vertices of the homotopy pushout above are $\cal F$-equivalent to certain representable functors through the respective $\cal F$-equivalences from $\cal F_2$. Testing by mapping into an arbitrary $\cal F$-local functor $W$, we find that the homotopy pushout above is $\cal F$-equivalent to the homotopy pushout of the respective representable functors.

The latter pushout is $\cal F$-equivalent to an representable functor $R_A$ through an $\cal F$-equivalence from $\cal F_3$.
\end{proof}

\section{Representability theorems}
We are ready now to prove the representability theorems.

\begin{theorem}\label{cohomological-rep}
Let $F\colon \cal S^{\op}\to \cal S$ be a small, homotopy functor converting coproducts to products, up to homotopy, and homotopy pushouts to homotopy pullbacks. Then there exists a fibrant simplicial set $Y$, such that $F(-)\simeq \cal S(-, Y)$. The value of $Y$ may be computed by substituting $\ast$ into $F$ and applying the fibrant replacement: $Y=\widehat{F(\ast)}$
\end{theorem}
\begin{proof}
We have proven so far that that the $Q$-localization constructed in \ref{Q-localization} is essentially the localization with respect to $\cal F$: every element of $\cal F$ is a $Q$-equivalence, hence $Q$-fibrant objects are $\cal F$-local and the inverse inclusion follows from Proposition~\ref{main_prop}, which says, in particular, that every $\cal F$-local object is also $Q$-fibrant, hence any $Q$-equivalence is also an $\cal F$-equivalence.

Given a small functor $F$ satisfying the conditions of the theorem, consider its fibrant replacement in the projective model structure $F\trivcofib \hat{F}$, then $\hat F$ is $\cal F$-local and therefore also $Q$-fibrant, hence the fibrant replacement of $\hat F$ in the $Q$-local model structure is a projective weak equivalence $F\simeq \hat F \we \cal S(-, \widehat{F(\ast)})$. Therefore it suffices to take $Y = \widehat{F(\ast)}$ to prove the first statement of the representability theorem. 

To construct an approximation by a cohomological functor for a functor $G$ consider the factorization of the map $G\to\ast$ into a trivial cofibration followed by a fibration in the $Q$-local model structure: $G\trivcofib \hat G\fibr \ast$. Then the map $\gamma\colon G\trivcofib \hat G$ is initial, up to homotopy, beneath maps of $G$ into other fibrant cohomological functors
\end{proof}

\begin{remark}\label{approximation by a cohomology functor}
Actually, we have proven a little bit more: for every functor $G\colon \cal S^{\op}\to \cal S$ there exists an approximation of $G$ by a universal, up to homotopy, cohomological functor, i.e., there exists a natural transformation $\gamma\colon G\to \hat G$, where $\hat G$ is cohomological, such that for every fibrant cohomological functor $H$, any map $G\to H$ factors through $\gamma$ and the factorization is unique up to simplicial homotopy.
\end{remark}

\begin{remark}\label{different}
There is a different, simpler, approach to the classification of cohomological functors, which also does not use the assumption that the functor is small: given a simplicial cohomological functor $G\colon \cal S^{\op}\to \cal S$, consider the natural map $q\colon G(X)\to \cal S(X,G(\ast))$ obtained by adjunction from the natural map $X=\cal S(\ast, X)\to \cal S(G(X),G(\ast))$, which exists, in turn, since $G$ is simplicial. The map $q$ is an equivalence if $X=\ast$, which gives a basis for induction on the cellular structure of $X$ similar to Proposition~\ref{main_prop}. This approach is simpler, and more general (works for all functors, not necessarily small), but it does not give the benefit of representing, cohomological functors as fibrant objects in a model category on small functors. We owe this remark to T.~Goodwillie.
\end{remark}

\begin{remark}
A similar representability result was obtained by J.F.~Jardine \cite{Jardine-rep}. His enriched representability theorem applies to fairly general model categories satisfying the conditions analogous to the definition of a well-generated triangulated category, but the conditions demanded from the functor in this work are much more restrictive then ours: commutation with arbitrary homotopy colimits. The fact that we restricted these conditions only to coproducts and homotopy pushout allows us to call it the enriched \emph{Brown} representability. Our method can be extended to other model categories as well, including those that do not satisfy the conditions of Jardine's theorem. In \cite{Dual-Brown} we prove a similar representability theorem in the dual category of spectra.
\end{remark}

Homological Brown representability for space-valued functors is essentially Goodwillie's classification of linear functors. We choose, however, to discuss the contravariant version of this theorem in our work (our result is related to Goodwillie's theorem in the same way as Adams' representability theorem \cite{Adams} related to G.W.~Whitehead's \cite{GWWhitehead} classification of generalized homology theories). Even though philosophically the two versions are the same, in order to obtain an implication between them, we would have to work out a stable analogue of our theorem and then use $S$-duality. We leave it to the interested reader.

\begin{definition}\label{def-homological}
Simplicial functor $F\colon \cal S^{\op}\to \cal S$ is called \emph{homological} if $F$ converts homotopy pushouts of  finite simplicial sets to homotopy pullbacks.
\end{definition}

\begin{example}
Any functor of the form $H_{X,Y}(-)=X\times{\cal S(-,Y)}$ is homological; we would like to distinguish homological functors of the form $H_{\ast, Y}$, hence the next definition.
\end{example}

\begin{definition}
A homological functor $F$ is \emph{reduced} if $F(\emptyset) \simeq \ast$.
\end{definition}

Similarly to Lemma~\ref{sphere} we have
\begin{lemma}\label{homology-of-sphere}
Let $F$ be a reduced homological functor, then for all $n\geq 0$ there is a weak equivalence $F(\partial \Delta^{n})\simeq \cal S(\partial \Delta^{n},Y)$, where $Y$ is a fibrant simplicial set weakly equivalent to $F(\ast)$.
\end{lemma}
\begin{proof}
The statement is proved by induction on $n$. For $n=0$ there is a weak equivalence $F(\partial \Delta^{0})=F(\emptyset)\simeq \ast = \cal S(\emptyset, \widehat{F(*)}) = \cal S(\partial \Delta^{0}, \widehat{F(*)})$. 

Suppose that the statement is true for $n$, then $\partial \Delta^{n+1} \simeq \Delta^n \coprod_{\partial\Delta^n} \Delta^n$, hence $F(\partial \Delta^{n+1}) \simeq \holim(F(\Delta^n)\to {F(\partial\Delta^n)} \leftarrow F(\Delta^n)$. Lemma~\ref{homotopy_pullback} implies that $F$ is a homotopy functor, hence  
\begin{align*}
F(\partial \Delta^{n+1}) &\simeq \holim(F(\ast)\to \cal S(\partial\Delta^n, \widehat{F(*})) \leftarrow F(\ast))  && \text{(inductive assumption)}\\
                         &\simeq \holim(\cal S(\ast,\widehat{F(\ast)})\to {\cal S(\partial\Delta^n, \widehat{F(*)})} \leftarrow \cal S(\ast,\widehat{F(\ast)})) && (\ast \text{ is a unit in }\cal S)\\
                         &\simeq \cal S(\hocolim(\ast\leftarrow \partial \Delta^n\to \ast), \widehat{F(*)}) 
                         \simeq \cal S(\partial \Delta^{n+1}, \widehat{F(*)})
\end{align*}
\end{proof}

\begin{theorem}\label{homological-rep}
Let $F$ be a reduced homological functor $F\colon \cal S^{\op}\to \cal S$, then for all finite simplicial sets $K\in S$ there is a weak equivalence $F(K)\simeq \cal S(K,Y)$, where $Y$ is a fibrant simplicial set weakly equivalent to $F(\ast)$. 
\end{theorem}
\begin{proof}
It is possible to prove this theorem along the lines of the proof of Theorem~\ref{cohomological-rep}, but the model categories appearing on the way are all combinatorial and the required localizations are all with respect to sets of maps, so the model theoretical part of this result is standard and not so interesting. Instead we chose to use the approach of Remark~\ref{different}.

Since $F$ is a simplicial functor, similarly to Remark~\ref{different} there is a natural map $F(X)\to \cal S(X,F(\ast))$, which is a weak equivalence if $X=\ast$.  This is the base for cellular induction.

Let $X$ be a finite simplicial set, i.e., there is a finite chain of inclusions $\emptyset = X_0\to X_1\ldots X_a\to X_{a+1}\to \ldots X_k = X$, so that $X_{a+1}$ is obtained from $X_a$ by attaching a cell:
\[
\xymatrix{
\partial\Delta^n \ar[r]\ar@{^(->}[d]& X_a\ar[d]\\
\Delta^n \ar[r] & X_{a+1}.
}
\]
Applying $F$ we obtain a homotopy pullback 
\[
\xymatrix{
F(\partial\Delta^n) & F(X_a)\ar[l]\\
F(\Delta^n) \ar[u] & F(X_{a+1})\ar[l]\ar[u].
}
\]

Assuming, by induction, that $F(X_a) = \cal S(F(X_a),\widehat{F(\ast)})$ and using Lemma~\ref{homology-of-sphere} we obtain:
\begin{align*}
F(X_{a+1}) &\simeq \holim(\cal S(\ast, \widehat{F(\ast)})\to \cal S(\partial\Delta^n, \widehat{F(*})) \leftarrow \cal S(F(X_a), \widehat{F(*)}))\\
                         &\simeq \cal S(\hocolim(\ast\leftarrow \partial \Delta^n\to X_a), \widehat{F(*)}) 
                         \simeq \cal S(X_{a+1}, \widehat{F(*)}).
\end{align*}
After $k$ steps we obtain $FX\simeq \cal S(X,\widehat{F(*)})$.
\end{proof}

\section{An example of a non-class-cofibrantly generated model category}\label{non-cofib}

The model of spaces on the category of small contravariant functors, which we constructed in Section~\ref{model-categories}, has a very nice property: every object in it is weakly equivalent to an $\aleph_0$-small object --- the representable functor. Our initial motivation for looking into this model category was to use this property in order to construct some homotopical localizations with respect to certain classes of maps, since the set-theoretical difficulties do not constitute an obstruction in our model. However, another difficulty came up and we could not overcome it so far: the localized model category on $\cal S^{\cal S^{\op}}$ is not class cofibrantly generated, hence the standard methods for constructing localizations are not applicable. On the other hand, this is the first example of a non-class-cofibrantly generated model category arising in the  topological context. Examples of model categories featuring similar properties, but taking origin in abstract category theory appeared in \cite{AHRT}. 

There are two slightly different versions of the definition of the class-cofibrantly generated model categories. The first one demands that the domains and the codomains of the generating (trivial) cofibrations are $\lambda$-presentable, and the second one in more general demanding only that the (co)domains are $\lambda$-small with respect to cofibrations. This confusion probably has its origin in the difference between the combinatorial model categories by J.~Smith and the cellular model categories by P.~Hirschhorn. For example, the projective model structure on $\cal S^{\cal S^{\op}}$ is class-cofibrantly generated of the first kind, while the equivariant model structure on the maps of spaces $\cal S^{[2]}$ is class-cofibrantly generated only of the second kind. The respective localizations of these model categories constructed in this paper are not class-cofibrantly generated

In order to see that our model category is not class-cofibrantly generated we formulate a simple 
\begin{proposition}\label{fib-closed}
Let $\cat M$ be a class-cofibrantly generated model category such that the domains and the codomains of the generating trivial cofibrations are $\lambda$-presentable for some cardinal $\lambda$. Then the fibrations are closed in the category $\mor{\cat M}$ under sequential $\lambda$-filtered colimits, in particular the fibrant objects are closed in $\cat M$ under sequential colimits. If the (co)domains of the generating trivial cofibrations are $\lambda$-small with respect to cofibrations only, then the same conclusion holds for sequential colimits with cofibrations as bonding maps. 
\end{proposition}
The proof is left to the reader.

If the localized model category on $\cal S^{\cal S^{\op}}$ would be class-cofibrantly generated, then the fibrant objects would be closed under sequential $\lambda$-filtered colimits by Proposition~\ref{fib-closed}. But it is easy to see that the representable functors are not closed under sequential colimits of any cardinality, hence the localized model category is not class-cofibrantly generated, at least by the first definition. 

Even more interesting example is the localization of the equivariant model category on $\cal S^{[2]}$. The fibrant objects (i.e., the diagrams equivariantly homotopy equivalent to  the orbits) are not closed under sequential colimits even if the bonding maps are cofibrations. Consider, for example, the following colimit: 
$\colim_{n<\omega} \left (\orbit{[n]}\right) = \left(\orbit{\aleph_0}\right)$, where $[n]=\coprod_{n}\ast$. It is quite surprising, but if we replace all the bonding maps by cofibrations, this colimit will be no longer equivalent to the orbit $\left(\orbit{\aleph_0}\right)$.
\[
\colim\left(
 \vcenter{
\xy
\def \leftshiftedbullet
{\save[]-<0.1cm,0cm>*{\bullet} \restore}
\xymatrix{
		&		&		&		&	&\save[]-<0.1cm,1.35cm>*{\bullet}\\
		&		&\save[]-<0.1cm,.7cm>*{\bullet}="t1" \restore &	&\save[]-<0.1cm,.7cm>*{\bullet} \ar@{-}@<.7cm>[r];[]-<.1cm,0cm>\restore & \save[]-<0.1cm,.7cm>*{\bullet}  \restore\\
{\bullet} \ar[d]^{\hphantom{!}}="c"     & {\bullet} \ar@{-}[r];[]+0_{}="a"  & \leftshiftedbullet & \bullet \ar@{-}[rr];[]+0&\leftshiftedbullet \ar@<-.1cm>[d]_{\hphantom{xx}}="e"^{\hphantom{x}}="e1" &\leftshiftedbullet\ar@{}@<1cm>[d]^{\hphantom{x}\ldots}="f"\\
{\bullet}                               & {\bullet} \ar@{-}[r];[]+0^{}="b"            & \leftshiftedbullet & \bullet \ar@{-}[rr];[]+0& \leftshiftedbullet& \save[]-<0.1cm,0cm>*{\bullet} \restore
\ar "a";"b"_{\hphantom{!}}="d"^{\hphantom{x}}="d1"
\ar@{^{(}->} "c";"d"^{}
\ar@{^{(}->} "d1";"e"
\ar@{^{(}->} "e1";"f"
}
\endxy
}\right)=
\vcenter{
\xy
\def \leftshiftedbullet
{\save[]-<0.1cm,0cm>*{\bullet} \restore}
\xymatrix{
&	&	&\save[]-<0.1cm,1.35cm>*{\bullet}& \ar@{-}@<1.35cm>[l]-<.1cm,0cm>_{\vdots }\\
 &	&\save[]-<0.1cm,.7cm>*{\bullet} \ar@{-}@<.7cm>[r];[]-<.1cm,0cm>\restore & \save[]-<0.1cm,.7cm>*{\bullet}  \restore& \ar@{-}@<.7cm>[l]-<.1cm,0cm>\\
 & \bullet \ar@{-}[rrr];[]+0_{}="g"&\leftshiftedbullet 
 & \leftshiftedbullet & \\ 
 & \bullet \ar@{-}[rrr];[]+0^{}="h"& \leftshiftedbullet & \save[]-<0.1cm,0cm>*{\bullet} \restore& , 
\ar "g";"h"
}
\endxy
}
\]
since if we try to map the orbit $\left(\orbit{\aleph_0}\right)$ into the last colimit, then such map must factor through one of the finite stages, and no map corresponds to the connected component of the identity map on $\orbit{\aleph_0}$ in the mapping space $\hom\left(\orbit{\aleph_0},\orbit{\aleph_0}\right)$. 

The same argument generalizes to sequential colimits of any cardinality, hence we can conclude that the localized model category on maps of spaces is not class-cofibrantly generated of the second kind.

\bibliographystyle{abbrv}
\bibliography{Xbib}

\begin{thebibliography}{10}

\bibitem{AHRT}
J.~Ad{\'a}mek, H.~Herrlich, J.~Rosick{\'y}, and W.~Tholen.
\newblock Weak factorization systems and topological functors.
\newblock {\em Appl. Categ. Structures}, 10(3):237--249, 2002.
\newblock Papers in honour of the seventieth birthday of Professor Heinrich
  Kleisli (Fribourg, 2000).

\bibitem{Adams}
J.~Adams.
\newblock A variant of {EH} {B}rown's represenability theorem.
\newblock {\em Topology}, 10:185--198, 1971.

\bibitem{BF:gamma}
A.~K. Bousfield and E.~M. Friedlander.
\newblock {Homotopy theory of $\Gamma$-spaces, spectra, and bisimplicial sets}.
\newblock In {\em {Geometric Applications of Homotopy Theory II}}, number 658
  in {Lecture Notes in Mathematics}. Springer, 1978.

\bibitem{Brown}
E.~H. Brown, Jr.
\newblock Cohomology theories.
\newblock {\em Ann. of Math. (2)}, 75:467--484, 1962.

\bibitem{CaCho}
C.~Casacuberta and B.~Chorny.
\newblock The orthogonal subcategory problem in homotopy theory.
\newblock In {\em An alpine anthology of homotopy theory}, volume 399 of {\em
  Contemp. Math.}, pages 41--53. Amer. Math. Soc., Providence, RI, 2006.

\bibitem{CSS}
C.~Casacuberta, D.~Scevenels, and J.~H. Smith.
\newblock Implications of the large-cardinal principles in homotopy theory.
\newblock Preprint, 1998.

\bibitem{PhDI}
B.~Chorny.
\newblock Localization with respect to a class of maps. {I}. {E}quivariant
  localization of diagrams of spaces.
\newblock {\em Israel J. Math.}, 147:93--139, 2005.

\bibitem{pro-spaces}
B.~Chorny.
\newblock A generalization of {Q}uillen's small object argument.
\newblock {\em Journal of Pure and Applied Algebra}, 204:568--583, 2006.

\bibitem{Dual-Brown}
B.~Chorny and G.~Biedermann.
\newblock Enriched brown representability for the dual category of spectra.
\newblock Preprint, 2011.

\bibitem{Chorny-Dwyer}
B.~Chorny and W.~G. Dwyer.
\newblock Homotopy theory of small diagrams over large categories.
\newblock {\em Forum Mathematicum}, 2007.
\newblock To appear.

\bibitem{DF}
E.~Dror~Farjoun.
\newblock Homotopy and homology of diagrams of spaces.
\newblock In {\em Algebraic topology (Seattle, Wash., 1985)}, Lecture Notes in
  Math. 1286, pages 93--134. Springer, Berlin, 1987.

\bibitem{Farjoun}
E.~Dror~Farjoun.
\newblock Homotopy theories for diagrams of spaces.
\newblock {\em Proc. Amer. Math. Soc. 101}, pages 181--189, 1987.

\bibitem{DZ}
E.~Dror~Farjoun and A.~Zabrodsky.
\newblock Homotopy equivalence between diagrams of spaces.
\newblock {\em J. Pure Appl. Algebra 41(2-3)}, pages 169--182, 1986.

\bibitem{Goo:calc2}
T.~G. Goodwillie.
\newblock Calculus. {II}. {A}nalytic functors.
\newblock {\em $K$-Theory}, 5(4):295--332, 1991/92.

\bibitem{Heller}
A.~Heller.
\newblock On the representability of homotopy functors.
\newblock {\em Jornal of the London Mathematical Society}, s2-23(3):551--562,
  1981.

\bibitem{Hirschhorn}
P.~S. Hirschhorn.
\newblock {\em Model categories and their localizations}, volume~99 of {\em
  Mathematical Surveys and Monographs}.
\newblock American Mathematical Society, Providence, RI, 2003.

\bibitem{Hovey}
M.~Hovey.
\newblock {\em Model categories}.
\newblock Mathematical Surveys and Monographs 63. American Mathematical
  Society, Providence, RI, 1999.

\bibitem{Jardine-rep}
J.~F. Jardine.
\newblock Representability theorems for presheaves of spectra.
\newblock {\em Journal of Pure and Applied Algebra}, 215(1):77--88, January
  2011.

\bibitem{Kelly}
G.~M. Kelly.
\newblock {\em Basic concepts of enriched category theory}, volume~64 of {\em
  London Mathematical Society Lecture Note Series}.
\newblock Cambridge University Press, Cambridge, 1982.

\bibitem{GWWhitehead}
G.~Whitehead.
\newblock Generalized homology theories.
\newblock {\em Trans. Amer. Math. Soc.}, pages 227--283, 1962.

\end{thebibliography}
\end{document}